\documentclass[12pt,a4paper]{article}
\usepackage[utf8]{inputenc}
\usepackage{polski}
\usepackage{latexsym,amsfonts,amssymb,amsmath,amsthm,amsxtra}
\usepackage{graphicx}

\begin{document}

\begin{center} {\large{\bf EMPIRICAL INTERPRETATION\\

\vspace{1mm}
OF THE PITMAN EFFICIENCY}}

\vspace{10mm}
Tadeusz Inglot\\

\vspace{5mm}
{\it Faculty of Pure and Applied Mathematics,\\ Wroc{\l}aw University of Science and Technology}
\end{center}

\vspace{4mm}
\begin{quotation}
\noindent{\bf Abstract.} We study an empirical interpretation of the Pitman efficiency in testing for uniformity in the two-parametric family of the beta distributions. We show that for contamination models the Pitman efficiency aproximates relative efficiency very well.

\vspace{2mm}
\noindent{\it Key words and phrases: Pitman efficiency, asymptotic relative efficiency, empirical relative efficiency, testing for uniformity, beta distribution.

\noindent{\it MSC Subject Classification:} 62G20, 62G10, 62F05.}

\end{quotation}

\vspace{3mm}
\noindent{\bf 1. Introduction.}

\vspace{2mm}
Usually, the Pitman efficiency is evaluated for subfamilies of a multidimiensional or nonparametric families of distributions indexed by a real parameter in some testing problem (e.g. Nikitin (1995), Lehmann and Romano (2008)). One can pose a question how these theoretical foundings are reflected in empirical behaviour of tests, under consideration.

The aim of the present note is to study an empirical interpretation of the Pitman efficiency in a particular example of testing for uniformity in the two-parametric family of the beta distributions.

\vspace{2mm}
\noindent{\bf 2. Pitman efficiency.}

\vspace{2mm}
We start with recalling a definition of the Pitman efficiency suitable for a future use. Let $\Gamma\subset \mathbb{R}^k,\;k\geqslant 1,$ be a nonempty set, ${\cal P}=\{P_{\gamma}:\gamma\in\Gamma\}$ a family of distributions on a measurable space $({\cal X},{\cal A})$ and $X_1,...,X_n$ a sample from a distribution $P\in {\cal P}$. Fix $\gamma_0\in\Gamma$. We test the null hypothesis $H_0:P=P_{\gamma_0}$ against $H_1:P\neq P_{\gamma_0}$. Suppose we want to compare two upper-tailed tests given by statistics $T_n,\,V_n$. For $0<\alpha<\beta<1$ and an alternative $P_{\gamma},\;\gamma\neq \gamma_0,$ let $N_T(\alpha,\beta,P_{\gamma})$ denote the minimal sample size such that for all $n\geqslant N_T(\alpha,\beta,P_{\gamma})$ the power of the test $T$ for the alternative $P_{\gamma}$ at the significance level $\alpha$ and for the sample size $n$ is not smaller than $\beta$. Similarly we define $N_V(\alpha,\beta,P_{\gamma})$ for the test $V$. The relative efficiency of $T$ with respect to $V$ is defined to be (cf. Nikitin (1995))

\vspace{2mm}
\hspace*{4.7cm}$\displaystyle {{\cal R}{\cal E}}_{TV}(\alpha,\beta,P_{\gamma})=\frac{N_V(\alpha,\beta,P_{\gamma})}{N_T(\alpha,\beta,P_{\gamma})}.$

\vspace{2mm}
Let $\gamma(s)\in\Gamma,\;s\in [0,1],$ be a continuous curve in $\Gamma$ such that $\gamma(0)=\gamma_0$. Assume that there exists a $\sigma$-finite measure $\lambda$ on $({\cal X},{\cal A})$ such that $P_{\gamma(s)}\ll\lambda$ for $s\in[0,1]$ and $H(P_{\gamma(s)},P_{\gamma_0})\to 0$ as $s\to 0,\,s>0$, where $H(P,Q)$  denotes the Hellinger distance between $P$ and $Q$. The family of distributions $\{P_{\gamma(s)}\}=\{P_{\gamma(s)}:s\in [0,1]\}$ we shall here call a path. The following definition is a version of those well-known from the literature (e.g. Nikitin (1995), Lehmann and Romano (2008)).

\vspace{2mm}
{\bf Definition.} Given $0<\alpha<\beta<1$ and a path $\{P_{\gamma(s)}\}$. If there exists the limit

\vspace{1.5mm}
\hspace*{2.7cm}$\displaystyle \lim_{s\to 0^+}{\cal R}{\cal E}_{TV}(\alpha,\beta,P_{\gamma(s)})=e_{TV}^P(\alpha,\beta,\{P_{\gamma(s)}\})\in [0,\infty]$,

\vspace{1.5mm}\noindent
then we call it the Pitman efficiency of the test $T$ with respect to the test $V$ for the path $\{P_{\gamma(s)}\}$.

\vspace{2mm}
Given a path $\{P_{\gamma(s)}\}$, consider the following assumption about a generic statistic $W_n$:

\vspace{1.5mm}\noindent
there exist scaling functions $\mu(s)\geqslant 0,$ $\sigma(s)>0,\;s\in[0,1],$ and a continuous distribution function $G(x)$ such that

\vspace{1.5mm}
\hspace*{3.9cm}$\displaystyle \lim_{n\to\infty}P_{\gamma_0}^n\left(\frac{W_n-\sqrt{n}\mu(0)}{\sigma(0)}\leqslant x\right)= G(x)$\hfill (1)

\vspace{1.5mm}\noindent
for all $x\in\mathbb{R}$ and for any sequence $s_n\to 0,\;s_n>0,$ we have

\vspace{2mm}
\hspace*{3.6cm}$\displaystyle \lim_{n\to\infty}P_{\gamma(s_n)}^n\left(\frac{W_n-\sqrt{n}\mu(s_n)}{\sigma(s_n)}\leqslant x\right)=G(x)$ \hfill (2)

\vspace{1mm}\noindent 
for all $x\in\mathbb{R}$.

\vspace{2mm}
Condition (2) is a little bit weaker than the uniform convergence in $s$ (cf. the condition (P1) in Serfling (1980) or the condition D in Noether (1955)).

Below we state a version of the Pitman theorem (cf. Lehmann and Romano (2008), Nikitin (1995)) in a form ready to apply in Sections 3 and 4. Its proof is rather known and is patterned on that of Kallenberg's Lemma (cf. Kallenberg (1983), Inglot (1999) Theorem 2.7 or \'Cmiel et al. (2019)). However, for the sake of completeness and for the convenience of the reader we provide it in the appendix. Usually (1) and (2) are fulfilled with $G(x)=\Phi(x)$ the standard normal distribution function. But in the theorem below and in its proof this fact is unimportant.\\

{\bf Theorem.} Suppose test statistics $T_n, V_n$ satisfy (1) and (2) for a path $\{P_{\gamma(s)}\}$ with the same distribution function $G(x)$ increasing on the set $\{x:0<G(x)<1\}$, functions $\sigma_T(s),\sigma_V(s)$ are right continuous at $s=0$, while $\mu_T(s), \mu_V(s)$ have nonnegative right-hand derivatives at $s=0$. Denote $c_T^P=(\mu'_T(0)/\sigma_T(0))^2$ and $c_V^P=(\mu'_V(0)/\sigma_V(0))^2$. If $\max\{c_T^P,c_V^P\}>0$ then there exists the Pitman efficiency of $T$ with respect to $V$ for the path $\{P_{\gamma(s)}\}$, does not depend on $\alpha$ and $\beta$ and equals

\vspace{0.1mm}
\hspace*{2cm}$\displaystyle e_{TV}^P(\{P_{\gamma(s)}\})=e_{TV}^P(\alpha,\beta,\{P_{\gamma(s)}\})=\left(\frac{\mu_T'(0)/\sigma_T(0)}{\mu'_V(0)/\sigma_V(0)}\right)^2=\frac{c_T^P}{c_V^P},$ \hfill (3)

\vspace{1mm}\noindent
where $c/0$ is understood as $\infty$.

\vspace{2mm}
\noindent{\bf 3. Testing problem and contamination paths.}

\vspace{2mm}
Consider a particular problem of testing for uniformity in the family of the beta distributions on the unit interval. 
So, now and in the sequel we restrict attantion to the following family of distributions

\vspace{1.5mm}
\hspace*{1.4cm}$\displaystyle {\cal P}=\{P_{\gamma}: P_{\gamma}=P_{(p,q,\varepsilon)}=(1-\varepsilon)P_{1\,1}+ \varepsilon P_{p\,q},\; \varepsilon\in[0,1], (p,q)\in \Theta\}$,

\vspace{1.5mm}\noindent
where $P_{p\,q}$ denotes the beta distribution on $[0,1]$ with parameters $p,q$, $\Theta=\{(p,q):\;p\geqslant q>0,\;\tau(p,q)\geqslant 0\}$ and $\tau(p,q)=2p^2-2pq-q^2+2p-q$. Let $\gamma_0=(1,1,0)$. Then $P_{\gamma_0}=P_{(1,1,0)}=P_{1\,1}$ is the uniform distribution. We test the simple null hypothesis $H_0:P=P_{1\,1}$ against $H_1:P\neq P_{1\,1},\;P\in{\cal P}$. Consider two upper-tailed tests given by the statistics $V_n=\sqrt{n}(\overline{X}-1/2)$ and $T_n=(\sum_{i=1}^n (X_i^2-1/3))/\sqrt{n}$.

Recall that

\vspace{2mm}
\hspace*{2.8cm}$\displaystyle E_{pq}X_1=\frac{p}{p+q},\;\;m_2=E_{pq}X_1^2=\frac{p(p+1)}{(p+q)(p+q+1)},$\hfill (4)

\vspace{1mm}
\hspace*{2.4cm}$\displaystyle m_4=E_{pq}X_1^4=\frac{p(p+1)(p+2)(p+3)}{(p+q)(p+q+1)(p+q+2)(p+q+3)}.$\hfill(5) 

\vspace{2mm}\noindent
Hence, for $P_{p\,q}$ with $\tau(p,q)<0$ we have $E_{pq}X_1^2=1/3+\tau(p,q)/[3(p+q)(p+q+1)]<1/3$ and for paths lying in the region $\tau(p,q)<0$ the statistic $T_n$ does not satisfy (1) and (2), as $\mu_T(s)<0$. Therefore we have required the restriction $\tau(p,q)\geqslant 0$ for the set $\Theta$.

\vspace{2mm}
Fix $(p,q)\in \Theta,\;(p,q)\neq(1,1),$ and let $P_{\gamma(s)}=(1-s)P_{1\,1}+sP_{p\,q},\;s\in[0,1]$. So, we have $\gamma(s)=(p,q,s)$ and the path $\{P_{\gamma(s)}\}$ links by a linear segment (in the space of distributions) $P_{1\,1}$ to $P_{p\,q}$. Hence $\{P_{\gamma(s)}\}$ forms a contamination family determined by a fixed distribution $P_{p\,q}$. Here we shall call it a linear path. Lyapunov's theorem and (4) imply that $V_n$ satisfies (1) and (2) with $G(x)=\Phi(x)$, $\mu_V(s)=(1-s)/2+sE_{pq}X_1-1/2=s(p-q)/2(p+q)$ and $\sigma_V^2(s)=(1-s)/3+ sm_2-(\mu_V(s)+1/2)^2$. The assumptions of Theorem are satisfied for this test and $\mu'_V(0)=(p-q)/2(p+q),\;\sigma_V(0)=1/\sqrt{12}$ and $c_V^P=3(p-q)^2/(p+q)^2$. Similarly, from (4), (5) and by Lyapunov's theorem it follows that $T_n$ satisfies (1) and (2) with $G(x)=\Phi(x)$, $\mu_T(s)=s\tau(p,q)/3(p+q)(p+q+1),\;\sigma_T^2(s)=(1-s)/5+s m_4-(\mu_T(s)+1/3)^2$. Hence the assumptions of Theorem are also satisfied for $T_n$ and $\mu'_T(0)=\tau(p,q)/3(p+q)(p+q+1)$, $\sigma_T^2(0)=4/45$, $c_T^P=5\tau^2(p,q)/4(p+q)^2(p+q+1)^2$. By (3) it follows that for $p>q$ with $\tau(p,q)>0$ the Pitman efficiency of $T$ with respect to $V$ for a linear path exists and equals

\vspace{2mm}
\hspace*{0.8cm}$\displaystyle {\cal E}_{TV}^P= {\cal E}_{TV}^P(P_{p\,q})=\frac{5(2p^2-2pq-q^2+2p-q)^2}{12(p-q)^2(p+q+1)^2}=\frac{5\tau^2(p,q)}{12(p-q)^2(p+q+1)^2}$.\hfill (6)

\vspace{2mm}\noindent
For $p=q<1$ the efficiency is equal to $\infty$ while for $p>q>1$ with $\tau(p,q)=0$ is equal to 0. Observe that for $p,q$ lying on the half-line given by the equation $p-2q+1=0,\; p>1,$ we have ${\cal E}_{TV}^P(P_{p\,q})=5/12$ while for $p,q$ on the segment $2p-q-1=0,\; 1/2<p\leqslant 1,$ we have ${\cal E}_{TV}^P(P_{p\,q})=5/3$. Consequently, for any path $\{P_{\gamma(s)}\}$contained in the above half-line (hence with $\varepsilon=1$) one has $e_{TV}^P(\{P_{\gamma(s)}\})=5/12$ while for any path contained in the above segment (i.e. with $\varepsilon=1$) one has $e_{TV}^P(\{P_{\gamma(s)}\})=5/3$.

{\small \begin{center} {\bf Table 1.} Empirical powers in \% (emp.p.) of the tests $V$ and $T$\\ and empirical relative efficiencies (${\cal ERE}$) for linear paths\\ determined by selected $P_{p\,q}$. $\alpha=0.05$, 100$\,$000 MC.

\vspace{1mm}
\begin{tabular}{|ll|c|r|rr|c|c|} \hline
\multicolumn{2}{|c|}{parameters}&&&\multicolumn{2}{|c|}{emp.p.}&&\\
\hspace*{1.5mm}$p$&\hspace*{1.5mm}$q$&$s(=\varepsilon)$&$n$\hspace*{1.2mm}&$V$&$T$\hspace*{0.6mm}&${\cal ERE}$&${\cal E}_{TV}^P$\\ \hline
5&4&1\hspace*{3.3mm}&80&56&0&&\\
&&&30000&100&0&$<0.003$&0\\ \hline
4&3.15&0.9   &80&{\bf 51}&02&&\\
&&            &3100&100&{\bf 51}&0.026&0.026\\ \hline
6&4&0.5       &100&{\bf 54}&15&&\\
&&            &462&99&{\bf 54}&0.216&0.220\\ \hline
3&1&0.2       &100&{\bf 54}&55&&\\ 
&&            &97&53&{\bf 54}&1.031&1.067\\ \hline
0.6667&0.5&0.5&200&{\bf 54}&78&&\\
&&            &100&35&{\bf 54}&2&2.074\\ \hline
0.55&0.5&0.9  &570&{\bf 54}&99&&\\
&&            &91&22&{\bf 54}&6.264&6.505\\ \hline
0.5&0.5&0.9   &36000&09&100&&\\
&&            &180&09&51&$>200$&$\infty$\\ \hline
\end{tabular}\end{center}}

\vspace{2mm}
Let us compare the theoretical formula (6) for linear paths with an empirical behaviour of both tests for several alternatives and the significance level 0.05. Since the distribution of the statistic $T_n$ has a little bit heavier right tail than the standard normal distribution we use an empirically obtained approximate formula $t_{\alpha n}=1.6445+0.1/(\sqrt{n+25}-5)$ for the critical value of the test $T$. It works well at least for sample sizes $n$ between 40 and 30000. Results are shown in Table 1. For each case the value of the parameter $s$ was chosen to keep an empirical power of the test $V$ close to 1/2 for a moderate sample size.

For the first case simulated powers of the test $T$ are close to 0.002 for both sample sizes and were marked in Table 1 as 0. It may be seen a very good proximity of empirical relative efficiencies and the Pitman efficiencies for all parameters $p$ and $q$ and all values of $s$, also far from 0.

\vspace{2mm}
\noindent{\bf 4. Other examples of paths.}

\vspace{2mm}
Take a path $\{P_{\gamma_1(s)}\}$ determined by $\gamma_1(s)=\gamma_1(s;p,q)=(1-s+ps,1-s-qs+2qs^2,1)$, where $p\geqslant q>0$ are fixed with $\tau(p,q)>0$ and $q<3+\sqrt{8}$. Here parameters of the beta distribution aproach (1,1) along a parabola. We do not contaminate them with $P_{1\,1}$. The condition for $q$ guarantees that $\gamma_1(s)$ is included in $\Gamma$.

Since densities of the beta distributions are continuous with respect to their parameters it follows by Lebesgue Bounded Convergence Theorem that the Hellinger distance $H(P_{\gamma_1(s)},P_{\gamma_0})$ tends to 0 when $s\to 0^+$ (parameters of the beta distibutions on $\gamma_1(s)$ are bounded away from 0 unifromly in $s$). From Lyapunov's theorem and (4) we get

\vspace{2mm}
\hspace*{0.3cm}$\displaystyle \mu_V(s)=\frac{1-s+ps}{2-2s+ps-qs+2qs^2}-\frac{1}{2}=\frac{(p+q)s-2qs^2}{2(2-2s+ps-qs+2qs^2)},\;\;s\in[0,1].$

\vspace{2mm}\noindent
Thus $\mu'_V(0)=(p+q)/4$. Since $\sigma^2_V(0)=1/12$ then $c_V^P=3(p+q)^2/4$. Similarly,

\vspace{2mm}
\hspace*{0.6cm}$\displaystyle \mu_T(s)=\frac{(1-s+ps)(2-s+ps)}{(2-2s+ps-qs+2qs^2)(3-2s+ps-qs+2qs^2)}-\frac{1}{3}$

\vspace{1mm}
\hspace*{1.7cm}$\displaystyle=\frac{(4p+5q+1)s}{3(2-2s+ps-qs+2qs^2)(3-2s+ps-qs+2qs^2)}+O(s^2),\;\;s\to 0^+.$

\vspace{2mm}\noindent
and hence $\mu'_T(0)=(4p+5q+1)/18$. As $\sigma^2_T(0)=4/45$ then $c_T^P=5(4p+5q+1)^2/144$. Finally by (3) of Theorem we obtain

\vspace{2mm}
\hspace*{4.5cm}$\displaystyle e_{TV}^P(\{P_{\gamma_1(s)}\})=\frac{5(4p+5q+1)^2}{108(p+q)^2}.$

\vspace{3mm}
For a simulation study take $p=6,q=4$. Then $(p_1(s),q_1(s))=(1+5s,1-5s+8s^2)$ form an arc of the parabola $q_1=0.32(p_1-1)^2-p_1+2$ and $e_{TV}^P(\{P_{\gamma_1(s;6,4)}\})=15/16=0.9375$. 

\begin{center}
{\small {\bf Table 2.} Empirical powers in \% (emp.p.) and empirical relative\\ efficiencies (${\cal ERE}$) of $T$ with respect to $V$ for selected alternatives\\ on the path $\{P_{\gamma_1(s;6,4)}\}$. $\alpha=0.05$, 100$\,$000 MC.

\vspace{1mm}
\begin{tabular}{|lll|c|r|rr|c|c|c|} \hline
\multicolumn{3}{|c|}{alternative}&&&\multicolumn{2}{|c|}{emp.p.}&&&\\
\hspace*{1.5mm}$s$&$p_1(s)$&$q_1(s)$&$\varepsilon$&$n$\hspace*{1.7mm}&$V$\hspace*{1mm}&$T$&${\cal ERE}$&${\cal E}_{TV}^P$&$e_{TV}^P$\\ \hline
0.5&3.5&0.5&0.1&170&{\bf 52}&62&&&\\
&&&&127&43&{\bf 52}&1.338&1.375&0.938\\ \hline
0.2&2&0.32&0.1&180&{\bf 52}&63&&&\\
&&&&131&42&{\bf 52}&1.374&1.420&0.938\\ \hline
0.1&1.5&0.58&0.2&130&{\bf 54}&60&&&\\
&&&&110&49&{\bf 54}&1.182&1.217&0.938\\ \hline
0.05&1.25&0.77&0.3&200&{\bf 54}&56&&&\\
&&&&189&52&{\bf 54}&1.058&1.083&0.938\\ \hline
0.02&1.10&0.9032&0.7&200&{\bf 52}&50&&&\\
&&&&206&52&{\bf 52}&0.971&0.996&0.938\\ \hline
0.01&1.05&0.9508&1\hspace*{2.6mm}&400&{\bf 53}&51&&&\\
&&&&419&55&{\bf 53}&0.955&0.967&0.938\\ \hline
\end{tabular}}\end{center}

\vspace{1mm}\noindent Choose several points on this path taking $s=0.5,\, 0.2,\, 0.1,$ $0.05,\, 0.02,\,0.01$. In Table 2 we show empirical powers and empirical relative efficiencies for selected alternatives on the path $\{P_{\gamma_1(s;6,4)}\}$ and values of ${\cal E}_{TV}^P$ for the linear paths corresponding to each selected point. Since both tests are very sensitive for all considered alternatives we contaminate them with $P_{1\,1}$ in the form $(1-\varepsilon)P_{1\,1}+\varepsilon P_{p_1(s)\,q_1(s)}$ with $\varepsilon$ chosen such that empirical powers of the test $V$ are close to 1/2 for moderate sample sizes.

For small $s=0.05$ and the sample size $200/e_{TV}^P(\{P_{\gamma_1(s;6,4)}\})=200\frac{16}{15}=213$ an empirical power of the test $T$ equals 0.58 thus significantly more than 0.54 for $V$.
From Table 2 it follows that for all cases empirical relative efficiencies are very close to ${\cal E}_{TV}^P$ (similarly as it was seen in Table 1). Even for quite small value $s=0.01$ (in this case $H(P_{\gamma_1(0.01)},P_{\gamma_0})\approx 0.045$), ${\cal E}_{TV}^P$ better approximates relative efficiency than $e_{TV}^P(\{P_{\gamma_1(s)}\})$.
Only for alternatives very close to $P_{1\,1}$, $e_{TV}^P(\{P_{\gamma_1(s)}\})$ can have a good empirical interpretation. But in this case both efficiencies are close each other. To see this by elementary calculations we get

\vspace{2mm}
\hspace*{3.4cm}$\displaystyle{\cal E}_{TV}^P(P_{p_1(s)\,q_1(s)})=\frac{5}{12}\left(\frac{4p+5q+1+o(s)}{3(p+q)+o(s)}\right)^2.$

\vspace{2mm}\noindent
Hence $\displaystyle\lim_{s\to 0^+} {\cal E}_{TV}^P(P_{p_1(s)\,q_1(s)})=e_{TV}^P(\{P_{\gamma_1(s)}\})$ which is nicely seen in two last columns of Table 2.

\vspace{2mm}
In simulations for $\{P_{\gamma_1(s)}\}$ we have had to contaminate alternatives with $P_{1\,1}$. So, to avoid a contamination consider two other paths $\{P_{\gamma_2(s)}\},\;\{P_{\gamma_3(s)}\}$ given by $\gamma_2(s)=(1+2s+s^2, 1+s+s^2,1)$ and $\gamma_3(s)=(1-s/2+s^2/2,1-2s/3,1)$. They link $P_{1\,1}$ to $P_{4\,3}$ and $P_{1\,1}$ to $P_{1\,1/3}$, respectively. Similarly as above, using (4) and (5) after simple calculations we get  for $\gamma_2(s)$: $\mu'_V(0)=1/4,\;\mu'_T(0)=1/6$ and $c_V^P=3/4,\;c_T^P=5/16,\; e_{TV}^P(\{P_{\gamma_2(s)}\})=5/12\approx 0.4167$ and for $\gamma_3(s)$: $\mu'_V(0)=1/24,\;\mu'_T(0)=2/27$ and $c_V^P=1/48,\;c_T^P=5/81,\;e_{TV}^P(\{P_{\gamma_3(s)}\})=80/27\approx 2.963$.

For a simulation study choose 5 alternatives on both paths taking $s=1,\,0.5,\,0.2,\,0.1,$ $0.05$. In Table 3 we show empirical powers and empirical related efficiencies for selected alternatives on the path $\{P_{\gamma_2(s)}\}$.
Now, a contamination is unnecessary which is marked in the table in the column denoted by $\varepsilon$. For better insight into the present example we add a column $H_2=H(P_{\gamma_2(s)},P_{\gamma_0})$ providing Hellinger distances  for considered alternatives.

\begin{center}
{\small {\bf Table 3.} Empirical powers in \% (emp.p.) and empirical related\\ efficiencies (${\cal ERE}$) of $T$ with respect to $V$ for selected alternatives\\ on the path $\{P_{\gamma_2(s)}\}$. $\alpha=0.05$, 100$\,$000 MC.

\vspace{1mm}
\begin{tabular}{|lll|l|c|r|rr|c|c|c|} \hline
\multicolumn{3}{|c|}{alternative}&&&&\multicolumn{2}{|c|}{emp.p.}&&&\\
\hspace*{1.5mm}$s$&$p_2(s)$&$q_2(s)$&$\varepsilon$&$H_2$&$n$\hspace*{1.7mm}&$V$\hspace*{1mm}&$T$\hspace*{0.6mm}&${\cal ERE}$&${\cal E}_{TV}^P$&$e_{TV}^P$\\ \hline
1&4&3&1&0.479&50&{\bf 57}&05&&&\\
&&&&&490&100&{\bf 57}&0.102&0.104&0.417\\ \hline
0.5&2.25&1.75&1&0.295&60&{\bf 52}&16&&&\\
&&&&&247&99&{\bf 52}&0.243&0.250&0.417\\ \hline
0.2&1.44&1.24&1&0.138&170&{\bf 52}&24&&&\\
&&&&&489&91&{\bf 52}&0.348&0.354&0.417\\ \hline
0.1&1.21&1.11&1&0.073&500&{\bf 51}&26&&&\\
&&&&&1305&87&{\bf 51}&0.383&0.387&0.417\\ \hline
0.05&1.1025&1.0525&1&0.038&1750&{\bf 51}&28&&&\\
&&&&&4375&85&{\bf 51}&0.400&0.402&0.417\\ \hline
\end{tabular}}\end{center}

Again for small $s=0.05$ and the sample size $1750/ e_{TV}^P(\{P_{\gamma_1(s)}\})=4200$ an empirical power of the test $T$ equals 0.50 thus less than 0.51 for $V$.

In Table 4 we present the results for $\{P_{\gamma_3(s)}\}$. In two last cases we have taken 10$\,$000 MC runs while in other cases 100$\,$000 MC runs. Again, we have added a column $H_3$ giving corresponding Hellinger distances.

For small $s=0.1$ and the sample size $7000/e_{TV}^P(\{P_{\gamma_1(s)}\})=2362$ an empirical power of the test $T$ equals 0.43 thus significantly less than 0.51 for $V$.
Results shown in Tables 3 and 4 confirm previous observations. Relative efficiencies are very well approximated by ${\cal E}_{TV}^P$ for all alternatives lying on two considered paths (now with $\varepsilon=1$) and rather poorly by $e_{TV}^P$ except very small values of $s$.

\begin{center}
{\small {\bf Table 4.} Empirical powers in \% (emp.p.) and empirical relative efficiencies (${\cal ERE}$)\\ of $T$ with respect to $V$ for selected alternatives on the path $\{P_{\gamma_3(s)}\}$. $\alpha=0.05$.

\vspace{1mm}
\begin{tabular}{|lll|l|c|r|rr|c|c|c|} \hline
\multicolumn{3}{|c|}{alternative}&&&&\multicolumn{2}{|c|}{emp.p.}&&&\\
\hspace*{1.5mm}$s$&\hspace*{1.5mm}$p_3(s)$&\hspace*{1.5mm}$q_3(s)$&\hspace*{1.4mm}$\varepsilon$&$H_3$&$n$\hspace*{1.7mm}&$V$\hspace*{1mm}&$T$\hspace*{0.6mm}&${\cal ERE}$&${\cal E}_{TV}^P$&$e_{TV}^P$\\ \hline
1&1&0.3333&0.2&0.161&100&{\bf 53}&64&&&\\
&&&&&73&44&{\bf 54}&1.370&1.437&2.963\\ \hline
0.5&0.875&0.6667&1&0.167&60&{\bf 56}&68&&&\\
&&&&&43&46&{\bf 56}&1.395&1.496&2.963\\ \hline
0.2&0.92&0.8667&1&0.055&1000&{\bf 50}&72&&&\\
&&&&&522&33&{\bf 50}&1.916&1.936&2.963\\ \hline
0.1&0.955&0.9333&1&0.027&7000&{\bf 51}&80&&&\\
&&&&&3080&30&{\bf 51}&2.273&2.295&2.963\\ \hline
0.05&0.97625&0.9667&1&0.013&38500&{\bf 51}&85&&&\\
&&&&&15010&28&{\bf 51}&2.565&2.58&2.963\\ \hline
\end{tabular}}\end{center}

Now, consider a path $\{P_{\gamma_4(s)}\}$ determined by $\gamma_4(s)=(1+s, 1+0.5 s-1.18 s^2,1)$. Both $\gamma_2(s)$ and $\gamma_4(s)$ are tangent at $s=0$ to $\ell(s)=(1+2s,1+s,1),\;s\in[0,1],$ where $(1+2s,1+s)$ form a segment of the line $p-2q+1=0$. Similar calculations, as previously, give $e_{TV}^P(\{P_{\gamma_4(s)}\})=5/12=e_{TV}^P(\{P_{\ell(s)}\})=e_{TV}^P(\{P_{\gamma_2(s)}\})$. 
To see how this fact can be observed empirically take three relatively close each other alternatives $P_{(1.2,1.1045,1)}$ lying on $\gamma_2(s)$ with $s\approx 0.095$, $P_{(1.2, 1.1,1)}$ lying on $\ell(s)$ with $s=0.1$ and $P_{(1.2,1.0528,1)}$ lying on $\gamma_4(s)$ with $s=0.2$. We have $H(P_{1.2\,1.1045},P_{1\,1})\approx 0.070,\;H(P_{1.2\,1.1},P_{1\,1})\approx 0.070,\;H(P_{1.2\,1.0528},P_{1\,1})\approx 0.077$ and $H(P_{1.2\,1.1045},P_{1.2\,1.1})\approx 0.002,\;H(P_{1.2\,1.0528},P_{1.2\,1.1})\approx 0.023$. In Table 5 we compare empirical relative efficiencies for these three alternatives.

\begin{center}
{\small {\bf Table 5.} Empirical powers in \% (emp.p.) and empirical relative\\ efficiencies (${\cal ERE}$) of $T$ with respect to $V$ for three alternatives\\ on above described paths. $\alpha=0.05$, 100$\,$000 MC.

\vspace{1mm}
\begin{tabular}{|c|ll|r|cc|c|c|c|} \hline
&\multicolumn{2}{|c|}{parametres}&&\multicolumn{2}{|c|}{emp.p.}&&&\\
path&\hspace*{2mm}$p$&\hspace*{5mm}$q$&$n$\hspace*{1.5mm}&$V$&$T$&${\cal ERE}$&${\cal E}_{TV}^P$&$e_{TV}^P$\\ \hline
$\gamma_2$&1.2&1.1045&550&{\bf 52}&27&&&\\
&&&1430&87&{\bf 52}&0.384&0.389&0.417\\ \hline
$\ell$&1.2&1.1&500&{\bf 52}&28&&&\\
&&&1215&85&{\bf 52}&0.412&0.417&0.417\\ \hline
$\gamma_4$&1.2&1.0528&220&{\bf 51}&37&&&\\
&&&350&69&{\bf 51}&0.629&0.637&0.417\\ \hline
\end{tabular}}\end{center}

Table 5 shows that for each case ${\cal E}_{TV}^P$ approximates relative efficiency very well while the Pitman efficiency $5/12$, common for these paths, says very little about them.

To confirm the above observation consider another path determined by $\gamma_5(s)=(1-s/2,1-2s/3,1)$ corresponding to the segment of the line $4p-3q-1=0$. Then $\gamma_3(s)$ and $\gamma_5(s)$ are tangent at $s=0$ and, similarly as previously, $e_{TV}^P(\{P_{\gamma_3(s)}\})=80/27=e_{TV}^P(\{P_{\gamma_5(s)}\})$. For an empirical comparison take two distributions close to $P_{1\,1}$ and close each other $P_{(0.95,0.9249,1)}$ lying on $\gamma_3(s)$ with $s\approx 0.1127$ and $P_{(0.95,0.9333,1)}$ lying on $\gamma_5(s)$ with $s=0.1$. Then $H(P_{0.95\,0.9249},P_{1\,1})\approx 0.030,\;H(P_{0.95\,0.9333},P_{1\,1})\approx 0.027$ and $H(P_{0.95\,0.9249},P_{0.95\,0.9333})\approx 0.005$. In Table 6 we compare empirical relative efficiencies for these alternatives.

\begin{center}
{\small {\bf Table 6.} Empirical powers in \% (emp.p.) and empirical relative\\ efficiencies (${\cal ERE}$) of $T$ with respect to $V$ for two alternatives\\ lying on above described paths. $\alpha=0.05$, 100$\,$000 MC.

\vspace{1mm}
\begin{tabular}{|c|ll|r|cc|c|c|c|} \hline
&\multicolumn{2}{|c|}{alternative}&&\multicolumn{2}{|c|}{emp.p.}&&&\\
path&\hspace*{2mm}$p$&\hspace*{5mm}$q$&$n$\hspace*{2.5mm}&$V$&$T$&${\cal ERE}$&${\cal E}_{TV}^P$&$e_{TV}^P$\\ \hline
$\gamma_3$&0.95&0.9249&5200&{\bf 51}&80&&&\\
&&&2350&30&{\bf 51}&2.213&2.241&2.963\\ \hline
$\gamma_5$&0.95&0.9333&12000&{\bf 51}&89&&&\\
&&&4120&26&{\bf 51}&2.913&2.918&2.963\\ \hline
\end{tabular}}\end{center}

Table 6 shows the same picture as Table 5 that $e_{TV}^P$, common for both alternatives, cannot distinguish and well approximate relative efficiencies which significantly differ while ${\cal E}_{TV}^P$ can.

The true (unknown) distribution $P\neq P_{1\,1}$ of the samle at hand (under $H_1$) lies on many paths. So, one cannot assign it a single number as its Pitman efficiency interpreted as a ratio of sample sizes guaranteing the same power of both tests. Even for alternatives quite close to $P_{1\,1}$ (in the sense of Hellinger distance) the Pitman efficiency can significantly differ from the relative efficiency as it depends on a shape of a considered path while does not depend on a value of the parameter $s$ corresponding to a selected alternative on it. This is the case when a path runs close to $P_{1\,1}$ for large values of $s$. For example, take $\gamma_6(s)=(1-s+1.1s^2,1-2s+2s^2,1)$ tangent to $m(s)=(1-0.25s,1-0.5s,1)$ at $s=0$ corresponding to the line $2p-q-1=0$ and $\gamma_7(s)=(1+0.2s-0.1s^2,1+0.1s-0.1s^2,1)$ tangent to $\ell(s)$ at $s=0$. Then we have $e_{TV}^P(\{P_{\gamma_6(s)}\})=5/3$ and $e_{TV}^P(\{P_{\gamma_7(s)}\})=5/12$. For $s=1$ we have $P_{\gamma_6(1)}=P_{\gamma_7(1)}=P_{(1.1,1,1)}$ with $H(P_{1.1\,1}, P_{1\,1})\approx 0.048$. For $P_{1.1\,1}$ an empirical power 0.50 is attained by the test $V$ for $n=400$ while by the test $T$ for $n=529$, resulting in an empirical relative efficiency 0.756. It is well approximated by ${\cal E}_{TV}^P(P_{1.1\,1})\approx 0.765$ but has nothing to do with the Pitman efficiencies 1.667 and 0.417, respectively.

\vspace{1mm}
Finally, consider a path given by $\gamma_8(s)=(2+s,1+s,s)$. So, for any $s$ we have 
contaminations of $P_{2+s\,1+2}$ with $P_{1\,1}$. The densities of $P_{2+s\,1+s}$ are uniformly bounded and 

\begin{center}
{\small {\bf Table 7.} Empirical powers in \% (emp.p.) and empirical relative\\ efficiencies (${\cal ERE}$) of $T$ with respect to $V$ for selected alternatives\\ on the path $\{P_{\gamma_8(s)}\}$. $\alpha=0.05$, 100$\,$000 MC.

\begin{tabular}{|lll|c|r|rr|c|c|c|} \hline
\multicolumn{3}{|c|}{alternative}&&&\multicolumn{2}{|c|}{emp.p.}&&&\\
$s=\varepsilon$&$p_8(s)$&$q_8(s)$&$H_8$&$n$\hspace*{1.2mm}&$V$\hspace*{1mm}&$T$&${\cal ERE}$&${\cal E}_{TV}^P$&$e_{TV}^P$\\ \hline
\hspace*{3mm}1&3&2&0.390&30&{\bf 64}&28&&&\\
&&&&77&98&{\bf 64}&0.390&0.417&0.938\\ \hline
\hspace*{1.5mm}0.8&2.8&1.8&0.262&30&{\bf 51}&27&&&\\
&&&&62&82&{\bf 51}&0.484&0.504&0.938\\ \hline
\hspace*{1.5mm}0.5&2.5&1.5&0.144&60&{\bf 52}&37&&&\\
&&&&96&70&{\bf 52}&0.632&0.651&0.938\\ \hline
\hspace*{1.5mm}0.2&2.2&1.2&0.054&300&{\bf 55}&48&&&\\
&&&&373&63&{\bf 55}&0.804&0.817&0.938\\ \hline
\hspace*{1.5mm}0.1&2.1&1.1&0.028&950&{\bf 51}&46&&&\\
&&&&1090&56&{\bf 51}&0.872&0.876&0.938\\ \hline
\hspace*{0.7mm}0.05&2.05&1.05&0.014&3600&{\bf 51}&48&&&\\
&&&&3990&55&{\bf 51}&0.902&0.907&0.938\\ \hline
\end{tabular}}\end{center}

\vspace{2mm}
\noindent uniformly bounded away from 0 with respect to $s\in [0,1/2]$ and converge to the density of $P_{2\,1}$ when $s\to 0^+$. The Pitman efficiency for this path equals $e_{TV}^P(\{P_{\gamma_8(s)}\})=15/16={\cal E}_{TV}^P(P_{2\,1})$ for the limiting distribution $P_{2\,1}$. In Table 7 it can be seen that again relative efficiencies are for each $s$ well approximated by ${\cal E}_{TV}^P(P_{2+s\,1+s})$ but not by $e_{TV}^P(\{P_{\gamma_8(s)}\})$ except very small $s$. It is easy to see that ${\cal E}_{TV}^P(P_{2+s\,1+s})\to 15/16=e_{TV}^P(\{P_{\gamma_8(s)}\})$ as $s\to 0^+$ which is seen in two last columns of Table 7. Additionally we give a column $H_8$ providing corresponding Hellinger distances.

\vspace{2mm}
\noindent{\bf Conclussion.}

\vspace{2mm}
The Pitman efficiency for linear paths approximates relative efficiency for any alternative on them very accurately while for arbitrary paths its empirical interpretation meets some essential difficulties.

\vspace*{5mm}
\noindent {\bf References}

\vspace{2mm}
\noindent Kallenberg, W. C. M. (1983). Intermediate efficiency, theory and examples, {\it Ann. Statist.}, 

{\bf 11} 1401-1420.

\vspace{1mm}
\noindent Inglot, T. (1999). Generalized intermediate efficiency of goodness of fit tests, {\it Math. Me-

thods Statist.}, {\bf 8} 487-509.

\vspace{1mm}
\noindent Inglot, T., Ledwina, T. and Ćmiel, B., (2019), Intermediate efficiency in nonparametric 

testing problems with an application to some weighted statistics, {\it ESAIM, Probability 

and Statistics}, {\bf 23} 697-738.

\vspace{1mm}
\noindent Lehmann, E. L. and Romano, J. P. (2008), {\it Testing Statistical Hypotheses}, Springer Texts 

in Statistics, Springer, New York, 3rd ed.

\vspace{1mm}
\noindent Nikitin, Y. (1995), {\it Asymptotic Efficiency of Nonparametric Tests}, Cambridge University 

Press, Cambridge.

\vspace{1mm}
\noindent Noether, G. E. (1955), On a theorem of Pitman, {\it Ann. Math. Statist.}, {\bf 26} 64-68.

\vspace{1mm}
\noindent Serfling, R., J. (1980), {\it Approximation Theorems of Mathematical Statistics}, Wiley, New 

York.\\

\noindent {\bf Appendix. Proof of Theorem.} 

\vspace{2mm}
Fix $0<\alpha<\beta<1$, a sequence $s_n\to 0^+$ and denote $t_{\alpha n}, v_{\alpha n}$ exact critical values of both compared tests. 

{\bf(i)} First we proof that $N_T(\alpha,\beta,P_{\gamma(s_n)})\to\infty$.

\noindent Let $p_0(x),p_s(x)$ denote denisties of $P_{\gamma_0}, P_{\gamma(s)}$ with respect to $\lambda$ and $\kappa_n=\kappa_{nT}=N_T(\alpha,\beta,P_{\gamma(s_n)})$. By $p_{n0}, p_{ns_n}$ we shall denote joint densities of the sample for distributions $P_{\gamma_0}, P_{\gamma(s_n)}$, respectively. For $A_n=\{T_{\kappa_n}> t_{\alpha\kappa_n}\}$ we have from Schwarz Inequality

\vspace{1.5mm}
\hspace*{0.95cm}$\displaystyle 0<\beta-\alpha\leqslant \int\limits_{A_n}(p_{\kappa_ns_n}-p_{\kappa_n 0})d\lambda^{\kappa_n}$

\vspace{1mm}
\hspace*{2.9cm}$\displaystyle=\int\limits_{A_n}(\sqrt{p_{\kappa_ns_n}}-\sqrt{p_{\kappa_n0}})^2d\lambda^{\kappa_n}+2\int\limits_{A_n}(\sqrt{p_{\kappa_ns_n}}-\sqrt{p_{\kappa_n0}})\sqrt{p_{\kappa_n0}}d\lambda^{\kappa_n}$

\vspace{1mm}
\hspace*{2.9cm}$\leqslant H^2(P_{\gamma(s_n)}^{\kappa_n},P_{\gamma_0}^{\kappa_n})+2\sqrt{\alpha}H(P_{\gamma(s_n)}^{\kappa_n},P_{\gamma_0}^{\kappa_n}).$\hfill (A1)

\vspace{1.5mm}\noindent
If $\kappa_n$ contains a bounded subsequence then, from the formula for Hellinger distance for product distributions i.e. $H^2(P_{\gamma(s_n)}^{\kappa_n},P_{\gamma_0}^{\kappa_n})=2-2(1-H^2(P_{\gamma(s_n)},P_{\gamma_0})/2)^{\kappa_n}$, the right hand side of (A1) would converge to 0 for this subsequence which contradicts the choice of $\alpha$ and $\beta$.  Therefore $\kappa_n\to\infty$.

{\bf (ii)} For the test $T$ we have $\displaystyle P_{\gamma_0}^{\kappa_n}\left(T_{\kappa_n}>t_{\alpha\kappa_n}\right)\leqslant\alpha\leqslant P_{\gamma_0}^{\kappa_n}\left(T_{\kappa_n}\geqslant t_{\alpha\kappa_n}\right).$ As $G(x)$ is increasing and continuous and $0<\alpha<1$ it follows from the assumption (1) of Theorem that 

\vspace{1.5mm}
\hspace*{4.5cm}$\displaystyle \frac{t_{\alpha\kappa_n}-\sqrt{\kappa_n}\mu_T(0)}{\sigma_T(0)}\to G^{-1}(1-\alpha)=z_{\alpha}.$

\vspace{2mm}\noindent 
Hence $t_{\alpha\kappa_n}=\sqrt{\kappa_n}\mu_T(0)+\sigma_T(0)z_{\alpha}+o(1)$.

{\bf (iii)} Denote $z_{\beta}=G^{-1}(1-\beta)$. We show that $s_n\sqrt{\kappa_n}$ converges to $(z_{\alpha}-z_{\beta})/\sqrt{c_T^P}$, if $c_T^P>0$, or to $\infty$, if $c_T^P=0$.

For $c_T^P>0$, suppose contrary that for some increasing sequence $k_n$ of integers it holds $s_{k_n}\sqrt{\kappa_{k_n}}\to C<(z_{\alpha}-z_{\beta})/\sqrt{c_T^P}$ and $\kappa_{k_n}$ is increasing. By the definition of $\kappa_n$ we have

\vspace{2mm}
\hspace*{0.7cm}$\displaystyle \beta\leqslant P_{\gamma(s_n)}^{\kappa_n}(T_{\kappa_n}> t_{\alpha\kappa_n})=P_{\gamma(s_n)}^{\kappa_n}(T_{\kappa_n}>\sqrt{\kappa_n}\mu_T(0)+\sigma_T(0)z_{\alpha}+o(1))$

\vspace{1mm}
\hspace*{10mm}$\displaystyle =P_{\gamma(s_n)}^{\kappa_n}\left(\frac{T_{\kappa_n}-\sqrt{\kappa_n}\mu_T(s_n)}{\sigma_T(s_n)}>\frac{(\mu_T(0)-\mu_T(s_n))\sqrt{\kappa_n}+\sigma_T(0)z_{\alpha}+o(1)}{\sigma_T(s_n)}\right)$.

\vspace{1mm}\noindent
Then for the subsequence $k_n$

\vspace{1mm}
\hspace*{-6mm}$\displaystyle \beta\leqslant P_{\gamma(s_{k_n})}^{\kappa_{k_n}}\left(\frac{T_{\kappa_{k_n}}-\sqrt{\kappa_{k_n}}\mu_T(s_{k_n})}{\sigma_T(s_{k_n})}>-\frac{\mu_T(s_{k_n})-\mu_T(0)}{s_{k_n}}\frac{s_{k_n}\sqrt{\kappa_{k_n}}}{\sigma_T(s_{k_n})}+\frac{\sigma_T(0)z_{\alpha}}{\sigma_T(s_{k_n})}+o(1)\right)$.\hfill (A2)

\vspace{1.5mm}\noindent
Define the sequence $\vartheta_n$ as follows: $\vartheta_{\kappa_{k_n}}=s_{k_n}$ and $\vartheta_j=s_j$ for $j\neq \kappa_{k_n}$. Then $\vartheta_n\to 0$ and from the assumption (2) of Theorem applied to the subsequence $\kappa_{k_n}$ of $\vartheta_n$ and from the existence of $\mu'_T(0)$, taking the limit in (A2) we get $\beta\leqslant 1-G(-C\sqrt{c_T^P}+z_{\alpha})$. This gives a contradiction since $z_{\alpha}-C\sqrt{c_T^P}>z_{\beta}$.

Similary supposing that for some increasing sequence $k_n$ of integers it holds $s_{k_n}\sqrt{\kappa_{k_n}}\to C'>(z_{\alpha}-z_{\beta})/\sqrt{c_T^P}$ we get form the minimality of $\kappa_n$ that

\vspace{2mm}
\hspace*{-6.5mm}$\displaystyle \beta > P_{\gamma(s_n)}^{\kappa_n-1}(T_{\kappa_n-1}> t_{\alpha\kappa_n-1})$

\vspace{1mm}
\hspace*{-6.5mm}$\displaystyle =P_{\gamma(s_n)}^{\kappa_n-1}\left(\frac{T_{\kappa_n-1}-\sqrt{\kappa_n-1}\mu_T(s_n)}{\sigma_T(s_n)}> -\frac{(\mu_T(s_n)-\mu_T(0))\sqrt{\kappa_n-1}-\sigma_T(0)z_{\alpha}+o(1)}{\sigma_T(s_n)}\right)$.\hfill (A3)

\vspace{2mm}\noindent
For the subsequence $k_n$ we have

\vspace{2mm}
\hspace*{-0.9cm}$\displaystyle \beta\geqslant P_{\gamma(s_{k_n})}^{\kappa_{k_n}-1}\left(\frac{T_{\kappa_{k_n}-1}-\sqrt{\kappa_{k_n}-1}\mu_T(s_{k_n})}{\sigma_T(s_{k_n})}>-\frac{\mu_T(s_{k_n})-\mu_T(0)}{s_{k_n}}\frac{s_{k_n}\sqrt{\kappa_{k_n}-1}}{\sigma_T(s_{k_n})}+z_{\alpha}+o(1)\right)$.\hfill (A4)

\vspace{2.5mm}\noindent
Define the sequence $\vartheta_n$ as follows: $\vartheta_{\kappa_{k_n}-1}=s_{k_n}$ and $\vartheta_j=s_j$ for $j\neq \kappa_{k_n}-1$. Then $\vartheta_n\to 0$ and from the assumption (2) of Theorem applied to the subsequence $\kappa_{k_n}-1$ after taking the limit in (A4) we get $\beta\geqslant 1-G(-C\sqrt{c_T^P}+z_{\alpha})$ which again gives a contradiction with $z_{\alpha}-C\sqrt{c_T^P}<z_{\beta}$.

If $c_T^P=0$, then suppose that for some increasing sequence $k_n$ of integers we have $s_n\sqrt{\kappa_n}\to C<\infty$. Repeating an analogous reasoning as in (A3), we obtain the inequality $\beta\leqslant 1-G(z_{\alpha})$ which also gives a contradiction.

{\bf (iv)} Similarly as in {\bf (iii)} we prove that $s_n\sqrt{\kappa_{nV}}\to (z_{\alpha}-z_{\beta})/\sqrt{c_V^P}$ if $c_V^P>0$ or $\infty$ if $c_V^P=0$. Hence, if $c_T^P,c_V^P>0$ then this and {\bf (iii)} give immediately

\vspace{2mm}
\hspace*{2cm}$\displaystyle \frac{N_V(\alpha,\beta,P_{\gamma(s_n)})}{N_T(\alpha,\beta,P_{\gamma(s_n)})}=\frac{s_n^2\kappa_{nV}}{s_n^2\kappa_{nT}}\longrightarrow \frac{c_T}{c_V}=\left(\frac{\mu'_T(0)/\sigma_T(0)}{\mu'_V(0)/\sigma_V(0)}\right)^2.$

\vspace{2mm}\noindent
If $c_T=0, c_V>0$ then $s_n\sqrt{\kappa_{nT}}\to\infty$ and simultaneously $s_n\sqrt{\kappa_{nV}}\to (z_{\alpha}-z_{\beta})/c_V$. Hence $\kappa_{nV}/\kappa_{nT}=s_n^2\kappa_{nV}/s_n^2\kappa_{nT}\to 0$. Similarly, if $c_T>0,c_V=0$ we get $\kappa_{nV}/\kappa_{nT}\to\infty$. This completes the proof of Theorem. \hfill$\Box$\\

\vspace{2mm}
\noindent Tadeusz Inglot\\
Faculty of Pure and Applied Mathematics\\
Wroc{\l}aw University of Science and Technology\\
Wybrze\.ze Wyspia\'nskiego 27, 50-370 Wroc{\l}aw, Poland.\\
E-mail: Tadeusz.Inglot@pwr.edu.pl

\end{document}